\documentclass[11pt]{article}
\usepackage{amsfonts}
\usepackage{latexsym,amssymb,amsfonts}
\usepackage{mathrsfs}
\usepackage{graphicx}
\usepackage{psfrag}
\usepackage{amsmath, amsbsy}
\usepackage{amsopn, amstext}
\usepackage{amsmath,multicol,amsopn}
\usepackage{latexsym,amssymb}
\def\sqr#1#2{{\vcenter{\vbox{\hrule height.#2pt
              \hbox{\vrule width.#2pt height#1pt \kern#1pt \vrule width.#2pt}
              \hrule height.#2pt}}}}
\def\signed #1{{\unskip\nobreak\hfil\penalty50
              \hskip2em\hbox{}\nobreak\hfil#1
              \parfillskip=0pt \finalhyphendemerits=0 \par}}
\def\endpf{\signed {$\sqr69$}}

\def\dbR{{\mathop{\rm l\negthinspace R}}}

\def\3n{\negthinspace \negthinspace \negthinspace }
\def\2n{\negthinspace \negthinspace }
\def\1n{\negthinspace }

\def\ds{\displaystyle}

\def\dbN{\mathbb{N}}

\def\dbR{\mathbb{R}}
\def\={\buildrel \triangle \over =}

\def\g{\gamma}

\def\e{\varepsilon}

\def\l{\lambda}

 \def\n{\nabla}

\def\f{\varphi}

\def\o{\omega}

\def\ns{\noalign{\ss} }

\def\pa{\partial}
\def\G{\Gamma}
\def\D{\Delta}

\def\O{\Omega}

\def\cA{{\cal A}}
\def\cB{{\cal B}}
\def\cC{{\cal C}}
\def\cD{{\cal D}}
\def\cE{{\cal E}}
\def\cF{{\cal F}}

\def\cJ{{\cal J}}

\def\cL{{\cal L}}
\def\cM{{\cal M}}

\def\cW{{\cal W}}
\def\cX{{\cal X}}

\def\no{\noindent}

\def\ms{\medskip}
\def\bs{\bigskip}
\def\q{\quad}
\def\qq{\qquad}
\def\hb{\hbox}

\def\mE{{\mathbb E}}

\def\max{\mathop{\rm max}}

\def\exp{\mathop{\rm exp}}

\def\pa{\partial}

\def\wt{\widetilde}
\def\cd{\cdot}

\def\as{\hbox{\rm a.s.{ }}}

\def\|{\Big |}
\def\({\Big (}
\def\){\Big )}
\def\[{\Big[}
\def\]{\Big]}
\def\be{\begin{equation}}
\def\bel{\begin{equation}\label}
\def\ee{\end{equation}}
\def\bt{\begin{theorem}}
\def\bcd{\begin{condition}}
\def\ecd{\end{condition}}
\def\et{\end{theorem}}
\def\bc{\begin{corollary}}
\def\ec{\end{corollary}}
\def\bde{\begin{definition}}
\def\ede{\end{definition}}
\def\bl{\begin{lemma}}
\def\el{\end{lemma}}
\def\bp{\begin{proposition}}
\def\ep{\end{proposition}}
\def\br{\begin{remark}}
\def\er{\end{remark}}
\def\ba{\begin{array}}
\def\ea{\end{array}}
\def\ed{\end{document}}
\def\ns{\noalign{\ms}}
\def\ds{\displaystyle}

\newcommand{\norm}[1]{\left\Vert#1\right\Vert}

\def\square#1{\vbox{\hrule\hbox{\vrule height#1%
     \kern#1\vrule}\hrule}}
\def\rectangle#1#2{\vbox{\hrule\hbox{\vrule height#1%
     \kern#2\vrule}\hrule}}

\font\tenbb=msbm10 \font\sevenbb=msbm7 \font\fivebb=msbm5

\newfam\bbfam
\scriptscriptfont\bbfam=\fivebb \textfont\bbfam=\tenbb
\scriptfont\bbfam=\sevenbb

\newtheorem{lemma}{Lemma}[section]
\newtheorem{remark}{Remark}[section]

\newtheorem{theorem}{Theorem}[section]
\newtheorem{corollary}{Corollary}[section]

\newtheorem{definition}{Definition}[section]
\newtheorem{proposition}{Proposition}[section]
\newtheorem{condition}{Condition}[section]

\makeatletter
   
   \@addtoreset{equation}{section}
\makeatother

\begin{document}
\title{\bf  Unique Continuation for Stochastic Heat Equations\thanks{This work is partially
ERC advanced grant 266907 (CPDENL) of the
7th Research Framework Programme (FP7), the
NSF of China under grant 11101070, the
project MTM2011-29306 of the Spanish
Science and Innovation Ministry,  the
Fundamental Research Funds for the Central
Universities in China under grants
ZYGX2012J115 and the General Fund Project
of Sichuan Provincial Department of
Education in China under the grant 124632.
\ms}
}
\author{ Qi L\"u\thanks{School of Mathematics, Sichuan University, Chengdu, 610064, China.
  {\small\it E-mail:} {\small\tt luqi59@163.com}.} \q and \q
  Zhongqi Yin\thanks{\textbf{Corresponding author}, School of Mathematics, Sichuan Normal University, Chengdu, 610068, China. {\small\it E-mail: }
  {\small\tt zhongqi.yin@yahoo.com}}}

\date{}

\maketitle

\begin{abstract}\no
We establish a unique continuation property  for stochastic heat
equations evolving in  a  domain $G\subset\dbR^n$($n\in\dbN$). Our
result shows that the value of the solution can be determined
uniquely by means of its value on an arbitrary open subdomain of $G$
at any given positive time constant. Further, when $G$ is convex and
bounded, we also give a quantitative version of the unique
continuation property. As applications, we get an observability
estimate for stochastic heat equations, an approximate result and a
null controllability result for a backward stochastic heat equation.
\end{abstract}

\bs

\no{\bf 2010 Mathematics Subject Classification}.  Primary 60H15,
93B05.

\bs

\no{\bf Key Words}. Stochastic heat equations, unique continuation
property, backward stochastic heat equations, approximate
controllability, null controllability. \ms

\section{Introduction}

\q\;We are concerned with the unique continuation for solutions to
stochastic heat equations and its application in approximate and null
controllability problems for backward stochastic heat equation. Let
us consider the following stochastic heat equation:
\begin{equation}{\label{system1}}
dy-\D ydt =  a  ydt  + b y  dB(t) \qq {\mbox { in }} G \times (0, T).
\end{equation}
Here, $T > 0$ is arbitrarily given and  $G
\subset \mathbb{R}^n (n\in \mathbb{N})$ is
a domain. Let $G_0\subset\subset G$ be a
subdomain.

We need to introduce some notations to be used in the  context.  By
$(\O, {\cal F}, \{{\cal F}_t\}_{t \geq 0}, P)$, we denote  a complete
filtered probability space, on which a standard one dimensional
 Brownian motion $\{ B(t) \}_{t\geq 0}$ is defined. Let $H$
be a Fr\'{e}chet space. We use the symbol $L^{2}_{\cal F}(0,T;H)$ to
stand for the Fr\'{e}chet space of all $H$-valued and $\{ {\cal F}_t
\}_{t\geq 0}$-adapted processes $X(\cdot)$ such that
$\mathbb{E}\big(\norm{X(\cdot)}^2_{L^2(0,T;H)}\big) < \infty$. By the
same manner, we adopt the symbols $L^{\infty}_{\cal F}(0,T;H)$ to
denote the Fr\'{e}chet space consisting of all $H$-valued
$\{\cF_t\}_{t\geq 0}$-adapted bounded processes and $L^{2}_{\cal
F}(\O;C([0,T];H))$ the Fr\'{e}chet space consisting of all $H$-valued
$\{\cF_t\}_{t\geq 0}$-adapted continuous processes $X(\cdot)$ such
that $\mE \big(\norm{X(\cdot)}^2_{C([0, T]; H)}\big)< \infty$. By
$L^2_{\cF_t}(\Omega; H)$, $0\leq t \leq T$,  we mean the Fr\'{e}chet
space consisting of  all $H$-valued $\cF_t$ measurable variables. All
the above spaces are endowed with the canonical quasi-norms.

We assume that the coefficients of the equation \eqref{system1}
satisfy the following condition:
$$
a\in L^{\infty}_{\cF}(0,T;L^{\infty}_{loc}(G)), \;\;  \;b\in
L^{\infty}_{\cF}(0,T; W^{1,\infty}_{loc}(G)).
$$

The definition to the solution of the equation \eqref{system1} is
given as follows.
\begin{definition}\label{def sol1}
We call $y\in L_{\cF}^2(\O;C([0,T];L^2_{loc}(G)))\cap
L_{\cF}^2(0,T;H^1_{loc}(G))$  a solution to the equation
\eqref{system1} if for any $t\in[0,T]$, any open set
$G'\subset\subset G$ and any $\f\in H_0^1(G')$, it holds that
\begin{equation}\label{def eq2}
\begin{array}{ll}\ds
\q\int_{G'}y(t,x)\f(x)dx-\int_{G'}y(0,x)\f(x)dx\\
\ns \ds = \int_0^t\int_{G'} \big[-\n y(s,x)\cd\n\f(x) +a(s,x)y(s,x)\f(x) \big]dx ds\\
\ns \ds\q +\int_0^t\int_{G'} b(s,x)y(s,x)\f(x) dx dB(s),\qq
P\mbox{-}\as
\end{array}
\end{equation}
\end{definition}

We have the following result.
\begin{theorem}\label{th-1}
For any given time $T_0\in (0,T]$, let $y$ be a  solution to the
equation \eqref{system1} such that $y(\cd,T_0)=0$ in $G_0$, $P$-a.s.
Then, $y(\cd,T_0)=0$ in $G$, $P$-a.s. Suppose, in addition, that $
a\in L^{\infty}_{\cF}(0,T;L^{\infty}(G))$, $b\in
L^{\infty}_{\cF}(0,T; W^{1,\infty}(G))$ and $y=0$ on $\partial G
\times (0,T)$, $P$-a.s. Then $y=0$ in $G\times (0,T)$, $P$-a.s.
\end{theorem}
\begin{remark}

If $b\in
L^\infty_\cF(0,T;W^{2,\infty}_{loc}(G))$
and $b_t\in
L^\infty_\cF(0,T;W^{2,\infty}_{loc}(G))$,
one can obtain Theorem \ref{th-1} in a
simple way.

Let
\begin{equation*}
z = e^{\ell} y, \quad \ell = - b(t,x) B(t).
\end{equation*}
In light of It\^o's formula, it follows that
\begin{equation*}
\begin{array}{ll}
\ds d z & = \ds y d e^{\ell} + e^{\ell} dy + dy d e^{\ell}\\
\ns & \ds = - y e^{\ell}[b_t B dt + b
dB(t)] + e^{\ell}[\D y dt + ay dt + by
dB(t)] - b^2 e^{\ell} y dt\\
\ns & \ds = (-b_t B z + e^{\ell}\D y  + a z  - b^2 z)dt\\
\ns & \ds = \D z dt - (b_t B - a + b^2 + B \D b ) z dt + 2 e^{\ell} B
\nabla b \cdot \nabla y
dt\\
\ns&\ds = \D z dt - (b_t B - a + b^2 + B \D b + |\n b|^2 B^2) z dt +
2  B \nabla b \cdot \nabla z dt.
\end{array}
\end{equation*}
Hence, we know that $z$ solves the following heat equation with
random coefficients
\begin{equation}\label{random-eq1}
z_t - \D z = 2 e^{\ell} B \nabla b \cdot \nabla z - (b_t B - a + b^2
+ B \D b + |\n b|^2B^2) z \q \mbox{ in } G\times (0,T).
\end{equation}
For a fixed $\o\in\O$, $z(\cd,\cd,\o)$ is a solution of a heat
equation.  From the unique continuation property of heat
equations(see \cite{Escauriaza1,phung1,Poon} for example), we have
that $z(\cd,T_0,\o)=0$ in $G$, provided that $z(\cd,T_0,\o)=0$ in
$G_0$. On the other hand, since $y(\cd,T_0,\cd)=0$ in $G_0$,
$P$-a.s., we have $z(\cd,T_0,\cd)=0$ in $G_0$, $P$-a.s. Thus, we
conclude that $z(\cd,T_0,\cd)=0$ in $G$, $P$-a.s., which leads to
that $y(\cd,T_0,\cd)=0$ in $G$, $P$-a.s. Further, if $y=0$ on $\pa
G\times (0,T)$, we get $z=0$ on $\pa G\times (0,T)$. This, together
with $z(\cd,T_0,\cd)=0$ in $G$, $P$-a.s.,  implies that $y=z=0$ in
$G\times (0,T)$.

Although the  strategy given in this remark is simple, it does not
work for our purpose in the present paper. Our method has its own
interest.

First, we can relax the regularity for the
coefficients. This is especially important
when one deals with semilinear equations.
Let us consider the following example.
\begin{equation}\label{semi-system1}
dw - \pa_{xx} wdt = w^m dB(t) \q \mbox{ in
} \;(0,L)\times (0,T).
\end{equation}
Here $L>0$ and $m\in\dbN$. From Theorem
\ref{th-1}, we can conclude that if
$$w_1,w_2\in C_\cF([0,T];H^2_{loc}(0,L))\cap
L^2_\cF(\O;C([0,T];H_{loc}^1(0,L)))$$
solves \eqref{semi-system1} and
$w_1(\cd,T_0)=w_2(\cd,T_0)$ in $G_0$,
$P$-a.s., then $w_1(\cd,T_0)=w_2(\cd,T_0)$
in $G$, $P$-a.s. Indeed, let $w_3=w_1-w_2$,
then $w_3$ solves
\begin{equation}\label{semi-system2}
dw_3 - \pa_{xx} w_3dt = \tilde b w_3 dB(t)
\q \mbox{ in } \;(0,L)\times (0,T),
\end{equation}
where $\tilde b = \sum_{j=0}^{m-1} w_1^j
w_2^{m-1-j}$. By Sobolev's embedding
theorem,  $\tilde b\in
L^\infty_\cF(0,T;W^{1,\infty}_{loc}(0,L))$.
Then,  by Theorem \ref{th-1},
$w_3(\cd,T_0)=0$ in $G$, $P$-a.s. Clearly,
$\tilde b$ will never be absolutely
continuous in the time variable $t$. Hence,
$\tilde b_t\notin
L^\infty_\cF(0,T;W^{2,\infty}_{loc}(0,L))$.

Second, if $G$ is convex, by our method, we can establish a
quantitative version of the unique continuation property. It seems
that this cannot be obtained by the strategy introduced above.
Indeed, as far as we know, to obtain an inequality like
\eqref{observability inequality} for the solution to the equation
\eqref{random-eq1}, we need the coefficients of $\n z$ and $z$ to be
essentially bounded with respect to $\o$. However, it is known that
the Brownian motion $B$ does not meet this condition.
\end{remark}

Further, if we put some more assumptions on $G$, $a$, $b$ and $y$, we
can get a better result than Theorem \ref{th-1}. More precisely, we
assume that
\begin{condition}\label{con1}
{\rm (1)} The domain $G$ is bounded and convex;
$$
(2)\;\; a\in L^{\infty}_{\cF}(0,T;L^{\infty}(G)), \;\; \;b\in
L^{\infty}_{\cF}(0,T; W^{1,\infty}(G)),  \quad y=0 \text{ on
}\partial G\times (0,T).
$$
\end{condition}

Under Condition \ref{con1}, the  solution to the equation
\eqref{system1} is now given in the following sense.

\begin{definition}\label{def sol}
We call $y\in
L_{\cF}^2(\O;C([0,T];L^2(G)))\cap
L_{\cF}^2(0,T;H^1_0(G))$ a solution to the
equation \eqref{system1} if for any
$t\in[0,T]$ and any $\f\in H_0^1(G)$, it
holds that
\begin{equation}\label{def eq2}
\begin{array}{ll}\ds
\q\int_{G}y(t,x)\f(x)dx-\int_{G}y(0,x)\f(x)dx\\
\ns \ds = \int_0^t\int_{G} \big[-\n y(s,x)\cd\n\f(x) +a(s,x)y(s,x)\f(x) \big]dx ds\\
\ns \ds\q +\int_0^t\int_{G} b(s,x)y(s,x)\f(x) dx dB(s),\qq
P\mbox{-}\as
\end{array}
\end{equation}
\end{definition}
Clearly, the $y$ satisfing Definition \ref{def sol} must satisfy
Definition \ref{def sol1}.

We have the following result.
\begin{theorem}\label{th-2}
Assume that Condition \ref{con1} holds. For any $T_0 \in (0, T]$,
there exist two constants $C > 0$ and $\delta \in (0, 1)$ such that
for any solution $y$ of \eqref{system1}, it holds that
\begin{equation}\label{th-2-eq}
\norm{y(T_0)}_{L^2(\Omega; L^2(G))} \leq C\norm{y(0)}_{L^2(\Omega;
L^2(G))}^{1 - \delta}\norm{y(T_0)}_{L^2(\Omega; L^2(G_0))}^{\delta}.
\end{equation}
\end{theorem}

\begin{remark}
Inequality \eqref{th-2-eq} is a quantitative version of the unique
continuation property for the solution $y$ with respect to $G_0
\times \{T_0\}$. Indeed, if $y$ solves the equation  \eqref{system1}
and $y(\cdot, T_0) = 0$ in $G_0$, $P$-a.s., then by \eqref{th-2-eq}
and the backward uniqueness of the stochastic heat equations (See
Lemma \ref{inv th1} given in the next section), it is clear that $y $
vanishes in $G$, $P$-a.s.

An obvious drawback of Theorem \ref{th-2}
is that we need  the convexity of $G$,
which is crucial in the proof. How to drop
it has its own independent interest.
\end{remark}

The research  of unique continuation for
solutions to partial differential equations
originated from the classical
Cauchy-Kovalevskaya theorem. It was studied
extensively in the literature. We refer the
readers to  \cite{Ho2,Zuily} and the rich
references therein  in this respect.
Besides its own interest in the partial
differential equation theory, it also plays
very important roles in both Inverse
Problems and Control Theory(see
\cite{LRS,Zu} for example). The classical
unique continuation property is of
qualitative nature. It guarantees that the
value of the solution in a given domain
$\cM_1$ can be uniquely determined  by that
of the solution in a suitable subdomain
$\cM_2$ of $\cM_1$. Once the unique
continuation property holds, a natural
question is whether one can find a way to
recover the solution in $\cM_1$ by the
values of the solution in $\cM_2$. It is
well known that the noncharacteristic
Cauchy problem is ill-posed, i.e., a small
error on the data in $\cM_2$ may cause
uncontrollable effects on the solution in
$\cM_1$ (see \cite{Hadamard} for example).
Therefore, it is important  to have
stability estimate for the solution. We
refer the readers to \cite{LRS} for an
introduction for this subject. In general,
such kind of estimates can be divided into
two classes:
\begin{enumerate}
  \item Observability estimate. The common form reads
  \begin{equation}
  \label{observability inequality}
    \norm{y}_{\cM_1} \leq \cC \norm{y}_{\cM_2},
  \end{equation}
  where $\norm{\cdot}_{\cM_i}, i=1,2$ denote some suitable norms for the
  restriction of the solution $y$ on $\cM_i, i=1,2$ respectively,
  and $\cC$ is a constant independent of $y$. We refer the readers to \cite{FI}
  and \cite{TZ} for the observability estimate for heat equations
  and stochastic heat equations, respectively.

  \item Quantitative version of unique continuation. It also appears
  in the form of \eqref{observability inequality} but with $\cC$
  depending on $y$. One can turn to \cite{Escauriaza1, Escauriaza2, Luqi,
  phung} for such kind of estimates for heat equations and its
  stochastic counterpart, respectively.
\end{enumerate}

Although there are numerous references addressing  to  unique
continuation properties for deterministic heat equations(see
\cite{Escauriaza1,Escauriaza2,Lin,phung,phung1,Poon,Vessella} for
example), very little is known for the stochastic counterpart and it
remains to be further understood.  As far as we know,
\cite{Luqi,Zhang1} are the only two publications in this field. The
result in \cite{Zhang1} shows that a solution to the stochastic heat
equation (without boundary condition) evolving in $G$ would vanish
almost surely, provided that it vanishes in ${G_0}\times (0,T)$,
$P$-a.s. In \cite{Luqi}, the authors proved that a solution to the
stochastic heat equation (with a partial homogeneous Dirichlet
boundary condition on arbitrary open subset $\G_0$ of $\pa G$)
evolving in $G$ vanishes almost surely, provided that its normal
derivative equals $0$ in ${\G_0}\times (0,T)$, $P$-a.s. Compared with
the result in Theorem \ref{th-1}, the results in \cite{Luqi,Zhang1}
do  not need the homogeneous Dirichlet boundary condition. However,
they have to utilize the information of the solution in the whole
time duration $(0,T)$. On the other hand, Theorem \ref{th-1} means
that the solution vanishes almost surely if it vanishes in
${G_0}\times \{T_0\}$. That is, we need only the information of the
solution to \eqref{system1} in $G_0$  at any fixed positive time
constant $T_0$ rather than the whole time duration.

As an application of Theorem \ref{th-1}, we give an approximate
controllability result for backward stochastic heat equations.  We
denote by $\{{\cal W}_t\}_{t \geq 0}$ the natural filtration
generated by $\{B(t)\}_{t \geq 0}$, which is augmented by all the
$P$-null sets. Let $H$ be a Banach space. We denote by $L^{2}_{\cal
W}(0,T;H)$  the Banach space of all $H$-valued and $\{ {\cal W}_t
\}_{t\geq 0}$-adapted processes $X(\cdot)$ with
$\mathbb{E}\big(\norm{X(\cdot)}^2_{L^2(0,T;H)}\big) < \infty$; by
$L^{\infty}_{\cal W}(0,T;H)$  the Banach space consisting of all
$H$-valued $\{\cW_t\}_{t\geq 0}$-adapted bounded processes and by
$L^{2}_{\cal W}(\O;C([0,T];H))$ the Banach space consisting of all
$H$-valued $\{\cW_t\}_{t\geq 0}$-adapted continuous processes
$X(\cdot)$ such that $\mE \big(\norm{X(\cdot)}^2_{C([0, T]; H)}\big)<
\infty$. All the above spaces are endowed with the canonical norms.

Assume that $\pa G$ is $C^2$. Consider the following controlled
backward stochastic heat equation
\begin{equation}\label{csystem2}
\left\{
\begin{array}{lll}\ds
dz + \D zdt =  a_1  zdt  + b_1Z dt + hdt + \chi_{E_1}\chi_{G_0}f dt +
ZdB(t) & {\mbox { in }} G \times (0, T),
 \\
\ns\ds z= 0 & \mbox{ on } \pa G \times (0, T),\\
\ns\ds z(T)=z_T &\mbox{ in } G.
\end{array}
\right.
\end{equation}
Here $E_1\subset [0,T]$ is a Lebesgue measurable subset with positive
measure, $z_T\in L^2(\O,\cW_T,P;L^2(G))$, $a_1\in
L^{\infty}_{\cW}(0,T;L^{\infty}(G))$, $b_1\in L^{\infty}_{\cW}(0,T;
W^{1,\infty}(G))$, $h \in L^2_\cW(0,T;L^2(G))$,   and $f \in
L^2_\cW(0,T;L^2(G))$ is the control.

Following the duality analysis in
\cite{Zhou},  one can show that
\eqref{csystem2} admits a unique solution
$$
(z,Z)\in \big[L^2_\cW(\O;C([0,T];L^2(G)))\cap
L^2_\cW(0,T;H_0^1(G))\big]\times L^2_\cW(0,T;L^2(G)).
$$
\begin{definition}\label{def app b heat}
System  \eqref{csystem2} is  {\it
approximately controllable} if for any
$z_T\in L^2(\O,\cW_T,P;L^2(G))$, any
$z_0\in L^2(G)$ and any $\e>0$, there
exists a control $f \in
L^2_{\cW}(0,T;L^2(G))$ such that the
solution to the system \eqref{csystem2}
with terminal state $z_T$ and control $f$
satisfying that
$\norm{z(0)-z_0}_{L^2(G)}\leq \varepsilon$.
\end{definition}

Making use of Theorem \ref{th-1}, we obtain the following result.
\begin{theorem}\label{th-3}
System \eqref{csystem2} is approximately controllable.
\end{theorem}

The approximate controllability for deterministic heat equations is a
classical topic in control theory and almost well understood now. We
refer the readers to \cite{T1,FPZ,FZ1,FZ2} and the rich references
therein for this topic. However, the approximate controllability
problems for forward and backward stochastic heat equations are quite
open. Some special cases in which the approximate controllability
problem for stochastic heat equations can be reduced to the same
problem for deterministic ones were studied in \cite{DMM,FGR,M1,R}.
In \cite{Luqi2}, the author shows that the null controllability of a
stochastic heat equation does not imply its approximate
controllability. As a direct consequence of the observability
estimate in \cite{TZ}, one can deduce the approximate controllability
of backward stochastic heat equations evolving in bounded domains.
Compared with the result in \cite{TZ}, in the present paper, the
control acts only on a measurable subset of $[0,T]$ rather than the
whole interval and the domain $G$ can be unbounded.

As another application of Theorem
\ref{th-2}, we have the following
observability inequality.

\begin{theorem}
\label{observability theorem} Assume that
Condition \ref{con1} holds. Let $E\subset
(0, T)$ be a subset with positive Lebesgue
measure and $G_0$ be a nonempty open subset
of $G$. Then any solution $y$ to
\eqref{system1} satisfies the estimate
\begin{equation}\label{ob-th-eq1}
  \norm{y(T)}^2_{L^2_{\cF_T}(\Omega; L^2(G))}\leq C \,\mE \int_E
  \int_{G_0}y^2(x,t)dxdt.
\end{equation}
\end{theorem}

Observability inequalities play important roles in the study of
controllability problems and state observation problems(see
\cite{Zu,Yamamoto2} for example). The observability estimate for heat
equations with lower order potentials depending both on $x$ and $t$
was first proved in \cite{FI} by a global Carleman estimate. In that
work, the integral with respect to $t$ is over $[0, T]$ rather than
only a measurable set. The inequality in the form of
\eqref{ob-th-eq1} was obtained in \cite{phung2} and was used to get
the null controllability and the Bang-bang principle for the time
optimal control of heat equations.

As an immediate consequence of Theorem \ref{observability theorem},
we can get the null controllability of  backward stochastic heat
equations. Let us first recall the following definition.

\begin{definition}\label{def null b heat}
System  \eqref{csystem2} is said to be {\it null controllable} if for
any $z_T\in L^2(\O,\cW_T,P;L^2(G))$, there exists a control $f \in
L^2_{\cW}(0,T;L^2(G))$ such that the solution of the system
\eqref{csystem2} with terminal state $z_T$ and control $f$ satisfying
that $z(0)=0$.
\end{definition}
We have the following result.
\begin{theorem}\label{null control}
Assume that $h=0$ in \eqref{csystem2}. Then
\eqref{csystem2} is null controllable,
provided that $G$ is bounded and convex.
\end{theorem}

The null controllability for heat equations
is  also a classical topic in the control
theory of partial differential equations.
There are a large of literatures about it
(see \cite{FG,FZ2,FI,Lebeau} and the
references therein). Similar to the case of
approximate controllability problems, the
study of null controllability problems for
forward and backward stochastic heat
equations are very uncomplete.   In
\cite{TZ}, the authors got the null
controllability of backward stochastic heat
equations evolving in bounded domains.
Compared with the result in \cite{TZ}, in
the present paper, the control acts only on
a measurable subset of $[0,T]$ rather than
the whole interval but the domain $G$
should be convex.

The rest of this paper is organized as follows: In Section
\ref{preliminaries}, we introduce some standing notations  and derive
some lemmas as preliminaries for the proof of Theorem \ref{th-1}.
Section \ref{proof of th1} is devoted to  the proof of Theorem
\ref{th-1}. Section \ref{proof of th2} is devoted to the proof of
Theorem \ref{th-2}. In Section \ref{proof of observability}, we prove
the observability result, i.e., Theorem \ref{observability theorem}.
Finally, Section \ref{proof of th3} is devoted to the proof of
Theorem \ref{th-3} and Theorem \ref{null control}.

\section{Some Preliminaries}
\label{preliminaries}

\q\;This section is devoted to some
preliminaries for the proof of Theorem
\ref{th-1}. In what follows, for simplicity
of notations, we adopt
$\vartheta(\cdot,x,x_0) =
\exp(-|x-x_0|^2/4(\cdot))$ and
$\norm{a}_{\cdot}, \norm{b}_{\cdot}$ for
$\norm{a}_{L^\infty_{\cF}(0,T;L^{\infty}(\cdot))}$,
$\norm{b}_{L^{\infty}_{\cF}(0,T;
W^{1,\infty}(\cdot))}$ respectively.. We
fix a point $x_0 \in {G_0}$. For $\lambda
> 0$, define
\begin{equation}\label{K}
K(x,t)\=(T - t + \lambda)^{-\frac{n}{2}}\vartheta(T - t+\lambda, x,
x_0),\quad (x, t) \in G\times [0, T].
\end{equation}
It is clear that
\begin{equation}
\label{for K}
\begin{cases}
K_t + \D K =0,\quad \nabla K = -\frac{x-x_0}{2(T-t+\lambda)}K,\\
\ns \ms \ds \D K = \frac{-n}{2(T-t+\lambda)}K + \frac{|x-
x_0|^2}{4(T -
t+\lambda)^2}K,\\
\ns \ms \ds K_{x_ix_j} = \frac{(x_i-x_{0i})(x_j - x_{0j})}{4(T -
t+\lambda)^2}K,\quad i \neq j.
\end{cases}
\end{equation}

\ms\ms

Let $\varphi$ be a $C^{\infty}$ function
with support $\wt G \subset G$ and $\Phi=
\varphi y$. Let $F = a\Phi - y \D \varphi -
2\nabla \varphi\cdot\nabla y$. It follows
that
\begin{equation}
  \label{Truncated sysstem}
d\Phi - \D \Phi dt = F dt + b\Phi dB.
\end{equation}

For $t\in [0,T]$, we put
\begin{equation}\label{HDN}
\left\{
\begin{array}{lll}
\ds H(t) =  \mE\int_G |\Phi (x,t)|^2 K(x,t) dx,\\
\ns\ds D(t) = \mE\int_G |\nabla \Phi(x,t)|^2 K(x,t) dx,\\
\ns\ds N(t) =  \frac{2D(t)}{H(t)} \;\mbox{ provided that }H(t) \neq
0.
\end{array}
\right.
\end{equation}

Throughout this section, we always work under the assumption
$H(\cdot)\neq 0$.

\begin{lemma}
  \label{lemma 1}
  For   the function $H(\cd)$ defined in \eqref{HDN}, involving the solution $y$ to the equation
  \eqref{system1}, it holds that
\begin{equation}
  \label{lemma 1 -estim}
H'(t) = -2D(t) + 2\mE \int_G \Phi F K dx + \mE \int_G b^2 \Phi^2 K
dx.
\end{equation}
\end{lemma}

{\it Proof}\,: Following It\^{o}'s formula and noticing that $\Phi$
has zero boundary condition, we have that
\begin{equation}
\label{for H}
\begin{array}{ll}
&\ds\quad H(t) - H(s) \\
\ns & = \ds 2\mE\int_s^t  \int \Phi d\Phi K dx + \mE\int_s^t\int_G
(d\Phi)^2 K dx + \mE\int_s^t\int_G \Phi^2 K_t dx d\tau\\
\ns \ds & = \ds 2\mE\int_s^t\!\int_G \Phi \big[\D  \Phi d\tau + F
d\tau + b\Phi dB(\tau)\big]K dx + \mE\int_s^t\!\int_G
(d\Phi)^2 K dx - \mE\int_s^t\int_G \Phi^2 \D K dx d\tau\\
\ns \ds & =\ds 2\mE\int_s^t\int_G \Phi\D \Phi K dxd\tau -
\mE\int_s^t\int_G \Phi^2 \D K
dx dt + 2\mE\int_s^t\int_G \Phi FK dx d\tau \\
\ns & \ds \quad + \mE\int_s^t\int_G (d\Phi)^2 K dx +
2\mE\int_s^t\int_G b\Phi^2 K dxdB(\tau)\\
\ns \ds & =\ds - 2\mE\int_s^t\int_G |\nabla\Phi|^2 K dx d\tau +
2\mE\int_s^t\int_G \Phi FK dx d\tau+ \mE\int_s^t\int_G (d\Phi)^2 K dx\\
\ns & \ds = -2\mE \int_s^t \int_G |\nabla \Phi|^2 K dx d\tau +
   2\mE \int_s^t \int_G \Phi F K dx d\tau + \mE \int_s^t \int_G b^2
   \Phi^2 K dx d\tau.
\end{array}
\end{equation}

As a result, it is easy to derive from \eqref{for H} that
\begin{equation*}
  H'(t) = -2 D(t) + 2\mE \int_G \Phi F K dx  + \mE \int_G b^2 \Phi^2
  K dx.
\end{equation*}
Thus, we complete the proof. \endpf

\begin{lemma}
  \label{lemma 2 }
  For $0\leq s < t \leq T$, it follows that
  \begin{equation}
    \label{lemma 2 estim}
\begin{array}
    {ll}
    \ds N(t) -N(s)\3n&\ds  \leq \int_s^t \(\frac{1}{T - \tau + \lambda} + 2\norm{b}^2_{\wt G}\)
    N(\tau)d\tau + 2\norm{b}^2_{\wt G}(t -s) + \int_s^t \frac{\mE \int_G F^2 K
    dx}{H}d\tau.
  \end{array}
  \end{equation}
\end{lemma}

{\em Proof}: First, we have that
\begin{equation}
\label{Lemm2-H}
\begin{array}{ll}
& \quad \ds D(t) - D(s)\3n\\
\ns & = \ds 2\mE \int_s^t\int_G \Phi d\Phi K dx + \mE \int_s^t\int_G
(d\Phi)^2 K dx  + \mE \int_s^t\int_G \Phi^2 K_t dx d\tau\\
\ns & = \ds 2\mE \int_s^t\int_G \Phi d\Phi K dx + \mE \int_s^t\int_G
(d\Phi)^2 K dx  - \mE \int_s^t\int_G \Phi^2 \D K dx d\tau\\
\ns & \ds = 2\mE \int_s^t\int_G \Phi d\Phi K dx  + \mE \int_s^t\int_G (d\Phi)^2 K dx\\
\ns & \ds \quad -\mE \int_s^t\int_G \nabla\cdot(\Phi^2 \nabla
K)dxd\tau + \mE \int_s^t\int_G
\nabla (\Phi^2)\cdot \nabla K dxd\tau\\
 \ns & =\ds 2\mE \int_s^t\int_G \Phi (\D \Phi + F)Kdx d\tau + 2\mE \int_s^t\int_G b \Phi^2 K dx dB(\tau)+ \mE \int_s^t\int_G
b^2\Phi^2 K dx d\tau \\
\ns &\ds \quad   - 2\mE \int_s^t\int_G \Phi \frac{x -
x_0}{2(T - \tau + \lambda)}\cdot\nabla \Phi K dx d\tau \\
\ns & =\ds 2\mE \int_s^t\int_G\Phi \(\D \Phi + \frac F 2 - \frac{x -
x_0}{2(T - \tau + \lambda)}\cdot \nabla
\Phi\) K dxd\tau \\
\ns & \ds \quad + \mE \int_s^t\int_G \Phi F K dx d\tau + \mE
\int_s^t\int_G b^2\Phi^2 K dxd\tau.
\end{array}
\end{equation}

\vspace{-0.2cm}

Next, following some straightforward
calculations, we have that
\begin{equation*}
  \begin{array}
  {ll}
  &\ds\quad  D(t) - D(s)\\
  \ns & = \ds 2\mE\int_s^t \int \nabla \Phi
  d\nabla
  \Phi K dx + \mE\int_s^t\int_G |d\nabla \Phi|^2 K dx  - \mE\int_s^t\int_G |\nabla \Phi|^2
  \D K dx d\tau\\
  \ns & \ds = -2 \mE\int_s^t\int_G d\Phi \D \Phi K dx + 2\mE\int_s^t\int_G d\Phi \frac{x - x_0}{2(T -\tau +\lambda)}\cdot \nabla \Phi K
  dx\\
  \ns & \ds \quad + \mE\int_s^t\!\int_G |d\nabla \Phi|^2 K dx \!- 2\mE\int_s^t\!\int_G \D
  \Phi \nabla \Phi \cdot\nabla K dx d\tau\! -\!\! 2\sum_{i=1}^n \mE\int_s^t\int_G
  \Phi_{x_i} \nabla \Phi\cdot\nabla K_{x_i} dxd\tau\\
  \ns & = \ds -2\mE\int_s^t\int_G d\Phi \D \Phi K dx + 2 \mE\int_s^t\int_G (d\Phi+\D \Phi d\tau) \frac{x - x_0}{2(T - \tau +
  \lambda)}\cdot \nabla \Phi K dx \\
  \ns \ds & \quad \ds + \mE\int_s^t\int_G |d\nabla \Phi|^2 K dx + \mE\int_s^t\frac 1{T - \tau + \lambda}\int_G |\nabla \Phi|^2 K
  dx d\tau \\
  \ns & \ds \quad - 2\mE\int_s^t\int_G \left(\frac{x - x_0}{2(T - \tau + \lambda)}\cdot\nabla
  \Phi\right)^2 K dx d\tau.
  \end{array}
  \end{equation*}

In the above equality, no boundary terms
appear due to the fact that $\Phi$ vanishes
in a neighborhood of the boundary of $G$.

Aided by the equation \eqref{Truncated sysstem} and regrouping the
right hand side of the above equality follows that
  \begin{equation}
  \label{Lemma-D1}
  \begin{array}
  {ll}
  & \ds\quad  D(t) -D(s) \\
  \ns & = \ds \mE\int_s^t \frac 1{T -\tau + \lambda}\int_G |\nabla \Phi|^2 K dx
  d\tau
  + \mE\int_s^t\int_G |d\nabla \Phi|^2 K dx \\
  \ns & \ds \quad - 2\[\mE\int_s^t\int_G \left(\frac{x - x_0}{2(T -\tau + \lambda)}\cdot\nabla \Phi\right)^2 K dx d\tau
   + \mE\int_s^t\int_G (\D \Phi + F) \D \Phi K dx d\tau \\
  \ns & \ds \quad - \mE\int_s^t\int_G (2\D \Phi + F)
  \frac{x - x_0}{2(T - \tau + \lambda)}\cdot\nabla \Phi K dx d\tau\]\\
  \ns & \ds \quad - 2\mE\int_s^t\int_G b\Phi \D \Phi K dx dB(\tau) + 2 \mE\int_s^t\int_G b\Phi
  \frac{x - x_0}{2(T -\tau + \lambda)}\cdot\nabla \Phi K dx dB(\tau)\\
  \ns & = \ds -2\mE\int_s^t\int_G \[ \D \Phi + \frac F
  2 - \frac{x - x_0}{2(T -\tau + \lambda)}\cdot \nabla \Phi\]^2 K dx d\tau + \frac 12 \mE\int_s^t\int_G F^2 K dxd\tau\\
  \ns & \ds \quad + \mE\int_s^t\frac 1{T - \tau + \lambda}\int_G |\nabla \Phi|^2 K dx d\tau +
  \mE\int_s^t\int_G |\nabla(b\Phi)|^2 K dxd\tau\\
  \ns & \ds \quad - 2\mE\int_s^t\int_G b\Phi \D \Phi K dx dB(\tau) + 2 \mE\int_s^t\int_G b\Phi
  \frac{x - x_0}{2(T -\tau + \lambda)}\cdot\nabla \Phi K dx dB(\tau).
  \end{array}
  \end{equation}

  On the other hand, we have that
  \begin{equation}
    \label{Lemma-D2}
    \begin{array}
      {ll}
      \ds D(t) \3n& =\ds \mE \int_G |\nabla \Phi|^2 K dx \\
      \ns & \ds = \mE \int_G \nabla\cdot (\Phi \nabla \Phi K )dx -
      \mE \int_G \Phi \D \Phi K dx - \mE \int_G \Phi \nabla
      \Phi\cdot\nabla K dx\\
      \ns & \ds = -\mE \int_G \Phi \D \Phi dx +\mE \int_G \Phi
      \nabla \Phi \cdot \frac{x - x_0}{2(T -t + \lambda)}K dx \\
      \ns & \ds = -\mE \int_G \Phi\left(\D \Phi +\frac F 2 - \frac{x - x_0}{2(T - t + \lambda)}\cdot \nabla
      \Phi\right)K dx + \frac 12 \mE\int_G F \Phi K dx.
    \end{array}
  \end{equation}

Based on \eqref{lemma 1 -estim}, \eqref{Lemma-D1} and
\eqref{Lemma-D2}, for any $0\leq s < t \leq T$, it follows that
\begin{equation*}
  \begin{array}
    {ll}
 & \quad \ds N(t) -N(s)\\
   \ns & = \ds 2\int_s^t \frac{HdD - D dH}{H^2}\\
   \ns & \ds = -\int_s^t \frac 4H\mE\int_G \[\D \Phi+\frac F 2 - \frac{x -x_0}{2(T - \tau+\lambda)}\cdot\nabla
   \Phi\]^2 K dx d\tau\\
   \ns & \ds \quad + \int_s^t \frac 4{H^2}\[\mE \int_G \Phi\(\D \Phi + \frac F2 - \frac{x-x_0}{2(T-\tau+\lambda)}\cdot \nabla \Phi\)K
   dx\]^2 d\tau\\
   \ns & \ds \quad - \int_s^t \frac 1{H^2}\(\mE \int_G F K dx\)^2
   d\tau - \int_s^t \frac{2D}{H}\mE \int_G b^2\Phi^2 K dx d\tau\\
   \ns & \ds \quad + \int_s^t \frac 2 H \mE \int_G |\nabla
   (b\Phi)|^2 K dx d\tau + \int_s^t \frac 1{T -\tau +\lambda}
   \frac{2D}{H} d\tau + \int_s^t \frac{\mE \int_G F^2 K dx}{H}d\tau.
  \end{array}
\end{equation*}

Applying Cauchy-Schwarz's inequality to the second term of the right
hand side of the above equality and noticing that $\Phi = \varphi y$
is supported in $\widetilde G$, we arrive at that
\begin{equation*}
  \begin{array}
    {ll}
    \ds N(t) - N(s) \3n& \ds \leq \int_s^t \frac {N(\tau)}{T -\tau
    +\lambda} d\tau + \int_s^t \frac 2 H \mE \int_G
    |\nabla(b\Phi)|^2 K dx d\tau + \int_s^t \frac{\mE \int_G F^2 K dx}{H}d\tau\\
    \ns & \ds \leq \int_s^t \(\frac 1{T -\tau + \lambda} +
    2\norm{b}^2_{\wt G}\)N(\tau)d\tau  + 2\norm{b}^2_{\wt G}(t -s) + \int_s^t
    \frac{\mE \int_G F^2 K dx}{H}d\tau.
  \end{array}
\end{equation*}

As a result, we complete the proof.\endpf

\ms

At last, we introduce the following backward uniqueness for
solutions to \eqref{system1}.

\begin{lemma}\label{inv th1}
Assume that Condition \ref{con1} holds. Let
$y$ be a solution to the equation
\eqref{system1}. If $y(T)=0$ in $G$,
$P$-a.s., then $y(t)=0$ in $G$, $P$-a.s.,
for all $t\in [0,T]$.
\end{lemma}

\begin{remark}
Lemma \ref{inv th1} was first proved in
\cite{Luqi3}(see \cite[Corollary
1.1]{Luqi3} for the details) for bounded
domain $G$. By a very similar but lengthy
argument, the same result can be obtained
for unbounded $G$. Hence, we omit the
detailed proof here.
\end{remark}

\section{Proof of Theorem \ref{th-1}}
\label{proof of th1}

This section is devoted to the proof of
Theorem \ref{th-1}. \ms

Let $B_{r_i}, i=1,2,3,4$ be balls with
center $x_0$ and radius $r_i$ respectively.
These balls have  the property $B_r =
B_{r_1}\subset G_0$, $
B_{r_j}\subset\subset
B_{r_{j+1}}\subset\subset G$($j=1,2,3$). We
also choose a special truncated function
$\varphi$ for the transformation $\Phi =
\varphi y$ such that $\varphi$ satisfies
the properties
\begin{equation*}
{\rm supp}\varphi \subset B_{r_4}\quad \text{and}\quad \varphi=1
\text{ in } B_{r_3}.
\end{equation*}

We borrow some ideas from \cite{phung,phung1}. In those papers, the
author defined functions like $N, H, D$ in our paper and discussed
the relations among them to obtain a quantitative unique
continuation for deterministic heat equations involving in a bounded
domain with Dirichlet boundary condition. But for the stochastic
settings, though the idea can be borrowed, the case is much
different for that there are random terms involved.

\ms

{\em Proof of Theorem \ref{th-1}}: Without loss
of generality and for the simplicity of the
notations, we only consider the case $T_0=T$
here. The general case can be handled similarly.

We first prove that if $y(x,T) = 0$, $P$-a.s. in $ B_{r_1}$, then it
also vanishes in $B_{r_2}$. We do this by contradiction argument. In
fact, if this claim does not hold, then $\mE \int_{B_{r_{2}}}
y^2(x,T)dx \neq 0$. Therefore, from the continuity of the solution
$y$ of \eqref{system1} with respect to $t$, there must be some small
positive constant $\varepsilon$, such that
\begin{equation*}
  \mE \int_{B_{r_2}} y^2(x,t)dx \neq 0\quad \text{for all } t\in [T - 2\varepsilon,
  T].
\end{equation*}
Let us consider the following inequality:
\begin{equation*}
  \int_G \|\nabla \[\Phi(x,T) \vartheta(2\lambda, x,x_0)\]\|^2 dx \geq
  0,
\end{equation*}
from which we can obtain some information about $N(T)$. By some
straightforward computations, one has that
\begin{equation}
  \label{N-13.1}
  \begin{array}
    {ll}
    &\q \ds \int_G \|\nabla \Phi(x,T) - \frac {x-x_0}{4\lambda}\Phi(x,T)\|^2
    \vartheta(\lambda, x,x_0)dx\\
    \ns & \ds \leq \int_G |\nabla \Phi(x,T)|^2 \vartheta(\lambda, x,x_0)dx + \int_G \frac{|x-x_0|^2}{16\lambda^2}
    \,\Phi^2(x,T) \vartheta(\lambda, x,x_0)dx\\
    \ns & \ds \quad - 2\int_G \frac{x - x_0}{4\lambda}\cdot \nabla
    \Phi(x,T)\,\Phi(x,T)\vartheta(\lambda, x,x_0)dx\\
    \ns &\ds \leq  \int_G |\nabla \Phi(x,T)|^2 \vartheta(\lambda, x,x_0)dx + \int_G \frac{|x-x_0|^2}{16\lambda^2}
    \,
    \Phi^2(x,T) \vartheta(\lambda, x,x_0)dx\\
    \ns & \ds \q +\frac{n}{4\lambda}\int_G
    \Phi ^2(x,T)\vartheta(\lambda, x,x_0)dx -\int_G
    \frac{|x-x_0|^2}{8\lambda^2}
    \Phi^2(x,T)\vartheta(\lambda, x,x_0)dx.
  \end{array}
\end{equation}
In the above process, we adopt integration
by parts in the last step. It is easy  to
find that
\begin{equation}
  \label{N-14.1}
  \begin{array}{ll}\ds
  \q\int_G \frac{|x - x_0|^2}{8\lambda} \Phi^2(x,T)\vartheta(\lambda, x,x_0)dx \\
  \ns\ds\leq \frac n2 \int_G \Phi^2(x,T)\vartheta(\lambda, x,x_0)dx + 2\lambda \int_G |\nabla \Phi(x,T)|^2 \vartheta(\lambda, x,x_0)dx.
  \end{array}
\end{equation}
Taking mathematical expectation from both sides of \eqref{N-14.1}
and noticing that
$$
\begin{array}{ll}\ds
&\ds \mE\int_G |\nabla \Phi(x,T)|^2 \vartheta(\lambda, x,x_0)dx \\
  \ns  = & \ds  \frac{\ds \mE\int_G |\nabla \Phi(x,T)|^2 \vartheta(\lambda, x,x_0)dx}
  {\ds \mE\int_G \Phi^2(x,T) \vartheta(\lambda, x,x_0)dx} \[\mE\int_G \Phi^2(x,T) \vartheta(\lambda, x,x_0)dx\]\\
  \ns \ds =  & \ds\frac 12  N(T) \mE\int_G \Phi^2(x,T) \vartheta(\lambda, x,x_0)dx,
\end{array}
$$
we have that
\begin{equation}
  \label{N-15.1}
  \begin{array}
    {ll}
    &\quad \ds \mE \int_G |x - x_0|^2 \Phi^2(x,T) \vartheta(\lambda, x,x_0)dx\\
    \ns   & \leq\ds  8\lambda \(\lambda N(T) + \frac n
    2\)\mE\int_G \Phi^2(x,T) \vartheta(\lambda, x,x_0)dx\\
  \ns & \leq \ds 8\lambda \(\lambda N(T) + \frac n
    2\)\bigg[\mE \int_{B_r} \Phi^2(x,T) \vartheta(\lambda, x,x_0)dx + \mE \int_{G\setminus B_r} \Phi^2(x,T_0) \vartheta(\lambda, x,x_0)dx\bigg]\\
  \ns & \leq \ds 8\lambda \(\lambda N(T) +   \frac n
    2\)\bigg[\mE   \int_{B_r}  \Phi^2(x,T_0) \vartheta(\lambda, x,x_0)dx \\
  \ns & \ds \quad  +  \frac 1{r^2} \mE   \int_{G\setminus B_r}  |x - x_0|^2 \Phi^2(x,T) \vartheta(\lambda, x,x_0)dx\bigg].
  \end{array}
\end{equation}

From this inequality, we immediately obtain that
\begin{equation}
  \label{N15.2}
\begin{array}
    {ll}
    &\quad \ds \[1- \frac{8\lambda}{r^2}\(\lambda N(T)+\frac{n}{2}\)\]\mE \int_G |x - x_0|^2 \Phi^2(x,T) \vartheta(\lambda, x,x_0)dx\\
    \ns   & \leq \ds 8\lambda \(\lambda N(T)+\frac{n}{2}\) \mE \int_{B_r} \Phi^2(x,T) \vartheta(\lambda, x,x_0)dx.
  \end{array}
\end{equation}
Now, \eqref{N15.2} shows that we should
give an estimate for $\lambda N(T)$.

\vspace{0.1cm}

Based on   Lemma \ref{lemma 2 } and Gronwall's
inequality, it follows for $t \in [T-
2\varepsilon, T]$ that
\begin{equation*}
  \begin{array}
    {ll}
\ds N(T) \3n& \ds \leq \[N(t) + 2T\norm{b}^2_{B_{r_4}} +\int_t^T
\frac{\mE\int_G F^2 K dx}{H}d\tau\] \exp\[\ds \int_t^T\(\frac 1{T
-\tau +\lambda} + 2\norm{b}^2_{B_{r_4}}\)d\tau\],
  \end{array}
\end{equation*}
which implies that
\begin{equation}
\label{Pr-1}
\begin{array}{ll}
\ds  \lambda N(T)\3n & \ds \leq (T +
\lambda)\exp\big(2T\norm{b}^2_{B_{r_4}}\big)
   N(t)\\
   \ns & \ds \quad + \;(T + \lambda)\exp\big(2T\norm{b}^2_{B_{r_4}}\big)\(2T\norm{b}^2_{B_{r_4}}+ \int_t^T \frac{\mE\int_G F^2 K
   dx}{H}d\tau\).
   \end{array}
\end{equation}
Integrating \eqref{Pr-1} with respect to
$t$ over $[T- 2\varepsilon, T -
\varepsilon]$, we find that
\begin{equation}
  \label{Pr-2}
  \begin{array}
  {ll}
\ds \varepsilon \lambda N(T) \3n& \ds \leq (T +
\lambda)\exp\big(2T\norm{b}^2_{B_{r_4}}\big)
  \int_{T- 2\varepsilon}^{T - \varepsilon} N(t)dt\\
   \ns & \ds \quad + \varepsilon(T + \lambda)\exp\big(2T\norm{b}^2_{B_{r_4}}\big)\(2T\norm{b}^2_{B_{r_4}}
   + \int_{T- 2\varepsilon}^T
\frac{\mE\int_G F^2 K
   dx}{H}d\tau\).
   \end{array}
\end{equation}
By  Lemma \ref{lemma 1}, we have that
\begin{equation}
  \label{Pr-3}
  \begin{array}
    {ll}
\ds \int_{T-2\varepsilon}^{T - \varepsilon}N(t)d t \leq
\ln\frac{H(T-2\varepsilon)}{H(T - \varepsilon)} +
2\int_{T-2\varepsilon}^{T - \varepsilon}\frac{\mE \int_G \Phi F K
dx}{H} dt + \varepsilon\norm{b}^2_{B_{r_4}}.
  \end{array}
\end{equation}
From \eqref{Pr-2} and \eqref{Pr-3}, we see
that
\begin{equation}
  \label{Pr-4}
  \begin{array}
    {ll}
    \ds \lambda N(T) & \ds \leq \frac{T +
    \lambda}{\varepsilon}\exp\big(2T\norm{b}^2_{B_{r_4}}\big)
    \[\ln\frac{H(T - 2\varepsilon)}{H(T - \varepsilon)} +
    2\int_{T-2\varepsilon}^{T- \varepsilon}\frac{\mE \int_G \Phi F K dx}{H} dt\\
    \ns & \ds \quad \varepsilon (1+ 2T)\norm{b}^2_{B_{r_4}} + \varepsilon \int_0^T
\frac{\mE\int_G F^2 K
   dx}{H}d\tau\]\\
   \ns & \leq \ds \frac{T +
    \lambda}{\varepsilon}\exp\big(2T\norm{b}^2_{B_{r_4}}\big)\[\ln\frac{H(T - 2\varepsilon)}{H(T - \varepsilon)}
     + \varepsilon + \varepsilon (1 + 2T)\norm{b}^2_{B_{r_4}}\\
    \ns & \ds \quad  + (\varepsilon + 1)\int_{T - 2\varepsilon}^T \frac{\mE\int_G F^2 K
   dx}{H}d t\].
  \end{array}
\end{equation}

For brevity of notations, we denote by
$\cA(\lambda)$ the right hand side of
\eqref{Pr-4}. We claim that the term
involving $F$ in $\cA(\lambda)$ are
uniformly bounded with respect to $\lambda
\in [0, 1]$ for fixed  $y$ solving
\eqref{system1}.


%
 Recalling that $\Phi =
\varphi y$ and $F = a\varphi y - 2\nabla
\varphi\cdot\nabla y - y\D \varphi$, we get that
\begin{equation*}
\begin{array}
{ll} \quad \ds \int_{T-2\e}^
T \frac{\mE \int_G F^2 K dx}{H}dt \\
\ns  \leq\ds 2T\norm{a}^2_{B_{r_4}}  + 8\int_{T-2\e}^T \frac{\mE
\int_G |\nabla \varphi|^2 |\nabla y|^2 K dx}{H}dt  +
2\int_{T-2\e}^T\frac{\mE \int_G |\D \varphi|^2 y^2 K}{H}dt \\
\ns \ds \leq 2T\norm{a}^2_{B_{r_4}} + 8\int_{T-2\e}^T \frac{\mE
\int_{G\setminus B_{r_3}} |\nabla \varphi|^2 |\nabla y|^2 K
dx}{\mE\int_{B_{r_2}} y^2 K dx}dt  + 2\int_{T-2\e}^T\frac{\mE
\int_{G\setminus B_{r_3}} |\D \varphi|^2 y^2 K }{\mE\int_{B_{r_2}}
y^2 K
dx}dt\\
\ns \ds \leq 2T\norm{a}^2_{B_{r_4}} + 8\max{|\nabla \f|^2}\exp\( -
\frac{r_3 - r_2}{4\lambda} \) \int_{T-2\e}^T
\frac{\mE\int_{G\setminus B_{r_3}} |\nabla y|^2
dx}{\mE\int_{B_{r_2}}y^2
dx}dt\\
\ns \ds \quad +  2\max{|\D
\f|^2}\exp\(-\frac{r_3 -
r_2}{4\lambda}\)\int_{T-2\e}^T\frac{\mE\int_{G\setminus
B_{r_3}} y^2 dx}{\mE\int_{B_{r_2}}y^2
dx}dt.
\end{array}
\end{equation*}
Note that $r_3 > r_2$, then for fixed $y$, we proved the claim that
the term involving $F$ is uniformly bounded with respect to
$\lambda\in[0, 1]$.


From the definition of $H(\cd)$, we see that
\begin{equation*}
\begin{array}
{ll}
 \ds \frac{H(T-2\varepsilon)}{H(T - \varepsilon)} \3n& =  \ds   \frac{\ds\mE\int_G \Phi^2(x,T - 2\varepsilon)K(x,T - 2\varepsilon)dx}
    {\ds\mE\int_G \Phi^2(x,T- \varepsilon)K(x, T - \varepsilon)dx}\\
    \ns  & \ds \leq \left(\frac{\lambda + \varepsilon}{\lambda + 2\varepsilon}\right)^{\frac
    n2}\frac{ \ds\mE\int_G
    \Phi^2(x,T - 2\varepsilon)\vartheta(\lambda+ 2\varepsilon, x,x_0)dx}
    {\ds\mE\int_G \Phi^2(x,T - \varepsilon)\vartheta(\lambda + \varepsilon, x,x_0)dx}.
 \end{array}
\end{equation*}

Denote the right hand side of the above inequality by
$\cE(\lambda)$. Then for fixed $\e$ and $\Phi$, it is clear that
$\cE(\lambda)$ is continuous in $(0,1]$ with respect to $\lambda$.
Then, we know that $\cE(\l)$ is uniformly bounded in $[0,1]$. So is
$\cA(\l)$.

Now returning to \eqref{N15.2}, noticing that $\lambda N(T) \leq
\cA(\lambda)$, we arrive at the following inequality:
\begin{equation}
\label{N-16.1}
\begin{array}{ll}\ds
  \q\[1- \frac{8\lambda}{r^2}\(\cA(\l)+\frac{n}{2}\) \]\mE \int_G |x - x_0|^2 \Phi^2(x,T) \vartheta(\lambda, x,x_0)dx\\
  \ns\ds \leq 8\lambda \(\cA(\l)+\frac{n}{2}\) \mE \int_{B_r} \Phi^2(x,T) \vartheta(\lambda, x,x_0)dx.
  \end{array}
\end{equation}
Since $\cA(\cd)$ is uniformly bounded in
$(0,1]$, we can choose  a $\l_1\in (0,1]$ such
that
$$
\[1- \frac{8\lambda_1}{r^2}\(\cA(\l_1)+\frac{n}{2}\) \]\geq \frac{1}{2}.
$$
This, together with \eqref{N-16.1}, implies that
\begin{equation*}
\mE \int_G |x - x_0|^2 \Phi^2(x,T) \vartheta(\lambda_1, x,x_0)dx \leq
r^2 \mE \int_{B_r} \Phi^2(x,T) \vartheta(\lambda_1, x,x_0)dx.
\end{equation*}
Returning to the definition of $\Phi = \varphi y$ and the special
properties of $\varphi$ and $r= r_1$, we arrive at that
\begin{equation}\label{5.14-eq1}
\mE \int_{B_{r_2}} |x - x_0|^2 y^2(x,T) \vartheta(\lambda_1, x,x_0)dx
\leq r^2_1 \mE \int_{B_{r_1}} y^2(x,T) \vartheta(\lambda_1, x,x_0)dx.
\end{equation}
From \eqref{5.14-eq1}, we know that we can
use the data in a small ball to estimate
that in a larger one. It follows that if
$y(x,T)$ vanishes in $B_{r_1}$, $P$-a.s.,
then so does in $B_{r_2}$.

Next, we show that if $y$, which is the
solution of \eqref{system1}, vanishes in a
given small ball, it may vanishes
everywhere in $G$. We prove this in the
following manner.

\ms

For brevity, we use $B$ in place of
$B_{r_1}$ in \eqref{5.14-eq1}. For any ball
$\widetilde B$ contained in $G$, we can
construct two sequences of balls
$\{S_i\}_{i=1}^m$ and $\{\widetilde
S_i\}_{i=1}^{m-1}$ containing in $G$ for
some finite natural number $m$, such that
\begin{equation*}
\left\{
\begin{array}
{ll}
  B\subset S_1, \mbox{ $B$ and $S_1$ have the same
  center,} \\
  \ns\ds  \widetilde S_i \subset \subset S_i\cap
  S_{i+1},\mbox{  $\widetilde S_i$ and $S_{i+1}$
  have the same center,
  }\; i=1,2,\cdots, m-1,\\
  \ns \ds  S_m\subset
  \widetilde B, \mbox{ $S_m$ and $\widetilde B$
  have the same center.}
\end{array}
\right.
\end{equation*}

By \eqref{5.14-eq1}, if $y(\cd,T)=0$ in
$B$, so does it in $S_1$. Then by the
selection of $\{S_i\}_{i=1}^m$ and
$\{\widetilde S_i\}_{i=1}^{m-1}$,
$y(\cd,T)=0$ in $\widetilde S_1$ follows.
So $y = 0$ in $S_2$, $\cdots$, and so on.
By induction, we conclude that if
$y(\cd,T)=0$ in $S_i$, so does it in
$\widetilde S_i \subset\subset S_i\cap
S_{i+1}, i=1,2,\cdots, m-1$. So $y(\cd,T)=
0$ in $S_{i+1}$ follows. we conclude that
$y=0$ in $S_m$. At last, noticing that $S_m
\subset \widetilde B$. Then use the same
argument to get \eqref{5.14-eq1}, we can
conclude $y(\cd,T)= 0$ in $\widetilde B$.

\vspace{0.1cm}

Following the arbitrariness of $\widetilde B$,
we conclude that if $y$, the solution of
\eqref{system1}, vanishes $P$-a.s. in a small
ball containing in $G$ at some fixed time $T$,
so it vanishes $P$-a.s in $G$ at time $T$.

\vspace{0.1cm}

Next, if $y=0$ on $\pa G\times (0,T)$, then
according to the backward uniqueness of
stochastic heat equations (See Lemma
\ref{inv th1}), we know $y=0$ in $G\times
(0,T)$, $P$-a.s. Recalling also that
$B_{r_1} = B_r \subset G_0$,  we complete
the proof.
\endpf

\section{Proof of Theorem \ref{th-2}}

\label{proof of th2}

In this section, we give the quantitative
unique continuation for the solution of the
equation \eqref{system1} subject to
Condition \ref{con1}. For convenience, we
adopt some notations.

For $K$ defined as in \eqref{K}, we
introduce that for $t\in [0,T]$,
\begin{equation}\label{HDN-convex}
\left\{
\begin{array}{lll}
\ds \widetilde H(t) =  \mE\int_G |y (x,t)|^2 K(x,t) dx,\\
\ns\ds \widetilde D(t) = \mE\int_G |\nabla y(x,t)|^2 K(x,t) dx,\\
\ns\ds \widetilde N(t) =  \frac{2\widetilde D(t)}{\widetilde H(t)}
\;\mbox{ provided that }\widetilde H(t) \neq 0.
\end{array}
\right.
\end{equation}

Similar to the process of establishing
Lemma \ref{lemma 1} and Lemma \ref{lemma 2
}, we have the following two results.

\begin{lemma}
  \label{lemma 3}
  For $\widetilde H$ defined in \eqref{HDN-convex}, it holds that
  \begin{equation}
    \label{H-convex}
    \widetilde H'(t) = - 2\widetilde D (t) + 2 \mE \int_G a y^2 K dx
    + \mE \int_G b^2 y^2 K dx.
  \end{equation}
\end{lemma}

{\em Proof}: The process is exactly
following that of proving Lemma \ref{lemma
1}.\endpf

\begin{lemma}
  \label{lemma 4}
 Assume that $\widetilde H \neq 0$.  For each $T > 0$ and $0 \leq s < t \leq T$, it follows that
  \begin{equation}
    \label{N-convex}
\begin{array}
    {ll}
    \ds \widetilde N(t) - \widetilde N(s)\3n&\ds  \leq \int_s^t \(\frac{1}{T - \tau + \lambda}+
    \norm{b}^2_{G}\)\widetilde N(\tau)d\tau + \int_s^t \( \norm{a}^2_{G}
    + 2\norm{b}^2_{G}\)d\tau.
  \end{array}
  \end{equation}
\end{lemma}

{\em Proof}\,: The process is mimicking that for Lemma \ref{lemma 2
}, so we omit some concrete calculations. The thing we should put
more attention to  is how to use the convex condition to deal with
the boundary term appearing in the proof.

First, we can arrive at that
\begin{equation}
\label{widetilde N-1}
  \begin{array}
    {ll}
    \widetilde H(t) - \widetilde H(s)& \!\!= \ds 2\mE \int_s^t \int_G y \(\D y +\frac 12 ay -\frac{x - x_0}{2(T - t + \lambda)}\cdot \nabla
   y  \)K dx d\tau + \mE \int_s^t\int_G a y^2 K
   dx d\tau \\
   \ns&\ds \q + \mE \int_s^t \int_G b^2 y^2 K
   dx d\tau.
  \end{array}
\end{equation}

Denote
\begin{equation}
\label{widetilde N-2}
  \cB = \nabla\cdot\(|\nabla y|^2 \nabla
    K\)- 2\nabla\cdot(\nabla y (\nabla y \cdot\nabla
    K)).
\end{equation}
We have that We have that
\begin{equation}
 \label{widetilde N-3}
  \begin{array}
    {ll}
   \ds \widetilde D(t) - \widetilde D(s) \3n & \ds = -2\mE \int_s^t\int_G \[\D y + \frac 12 ay - \frac{x-x_0}{2(T - t+\lambda)}\cdot\nabla
    y\]^2 K dx d\tau
    \\
    \ns & \ds \q + \frac 12 \mE\int_s^t \int_G a^2 y^2 K dxd\tau -\mE \int_s^t\int_G \cB dxdt\\
    \ns & \ds\q  + \frac{1}{T - t + \lambda}\int_G |\nabla y|^2 K
    dxdt + \int_G |\nabla(by)|^2K  dxdt.
  \end{array}
\end{equation}
We also have that
\begin{equation}\label{widetilde N-4}
\begin{array}{ll}\ds
\widetilde D(t) = -\mE\int_G y\[ \D y+ \frac{1}{2}ay  -
\frac{x-x_0}{2(T-t+\l)}\cdot\n y \]K dx + \frac{1}{2}\mE\int_G ay^2
K dx.
\end{array}
\end{equation}
According to \eqref{widetilde
N-1}--\eqref{widetilde N-4}, we can find
that for any $0\leq s < t \leq T$,
\begin{equation}
  \label{widetilde N}
  \begin{array}
    {ll}
   & \quad \ds \widetilde N(t) -\widetilde N(s)\\
   \ns & \ds =-\int_s^t \frac 4{\widetilde H}\mE\int_G \[\D y + \frac 12 ay - \frac{x - x_0}{2(T-\tau+\lambda)}\cdot \nabla y
    \]^2 K dx d\tau \\
    \ns & \ds \q
    -\int_s^t \frac 2{\widetilde H} \mE\int_G \cB dx d\tau  + \int_s^t \frac{1}{T - \tau + \lambda}\frac{2\widetilde D}{\widetilde H}d\tau
     +  \int_s^t\frac 2{\widetilde H} \mE\int_G |\nabla(by)|^2
    dxd\tau\\
    & \ds\q +  \int_s^t \frac{4}{\widetilde H^2}\,\[\mE\int_G  y \(\D y + \frac 12 ay -\frac{x - x_0}{2(T - \tau + \lambda)}
    \cdot \nabla
   y dt \)K dx \]^2 d\tau\\
   \ns & \ds \q - \int_s^t \frac 1{\widetilde H^2}\(\mE\int_G ay^2K dx\)^2d\tau - \int_s^t\frac{2\widetilde D}{\widetilde H^2} \mE\int_G b^2 y^2 K
   dx dt + \int_s^t \frac 1{\widetilde H} \mE\int_G a^2y^2 K dx d\tau.
  \end{array}
\end{equation}

For $G$ being convex, then $(x - x_0)\cdot
\nu \geq 0$ for each $x \in \pa G$ with
$\nu$ the out unit normal vector at $x$.
Also, noting that $y$, which solves the
equation \eqref{system1}, vanishes on the
boundary $\pa G$, one finds that
\begin{equation}
\label{N-boundary}
\begin{array}
{ll}
  \ds \int_s^t \mE \int_G \cB dxd\tau \3n& \ds = \int_s^t \mE\int_G \[\nabla \cdot(|\nabla y|^2\nabla K) - 2\nabla\cdot(\nabla y (\nabla y \cdot\nabla
    K))\] dxd\tau\\
    \ns & \ds =\int_s^t \mE\int_{\partial G}\[(|\nabla y|^2\nabla K) - 2(\nabla y (\nabla y \cdot\nabla
    K))\]\cdot \nu\;dS d\tau\\
    \ns & \ds =\!\int_s^t \!\mE \int_{\partial G}\[\!-|\nabla y|^2 \frac{(x - x_0)\cdot\nu}{2(T\!- \!\tau\! + \!\lambda)}
     + 2  \frac{x - x_0}{2(T\! - \!\tau \!+\! \lambda)}\cdot \nabla y (\nabla y \cdot
     \nu)\]K dS d\tau\\
     \ns & \ds = \int_s^t \mE \int_{\partial G}\[-|\nabla y|^2 \frac{(x - x_0)\cdot\nu}{2(T -\tau + \lambda)}
     + 2  |\nabla y|^2 \frac{(x - x_0)\cdot\nu}{2(T -\tau + \lambda)}\]K
     dS
     d\tau\\
     \ns & \ds \geq 0,
\end{array}
\end{equation}
where we use the symbol $dS$ to represent the boundary measure.

Combining \eqref{widetilde N} with
\eqref{N-boundary},  it follows that
\begin{equation}
   \label{N-convex}
  \begin{array}
    {ll}
 \ds \widetilde N(t) - \widetilde N(s)\3n & \ds \leq
       \int_s^t \frac{1}{T - \tau + \lambda}\frac{2\widetilde D}{\widetilde H}d\tau
     +  \int_s^t\frac 2{\widetilde H} \mE\int_G |\nabla(by)|^2
    dxd\tau + \int_s^t \frac 1{\widetilde H} \mE\int_G a^2y^2 K dx d\tau\\
    \ns & \ds  \leq \int_s^t \(\frac{1}{T - \tau + \lambda}+
    \norm{b}^2_{G}\)\widetilde N(\tau)d\tau + \int_s^t \( \norm{a}^2_{G}
    + 2\norm{b}^2_{G}\)d\tau.
  \end{array}
\end{equation}
 \endpf

\begin{remark}
For Theorem \ref{th-2}, we assume the convexity
of the domain $G$. We use this assumption  to
deal with the boundary terms as shown in
\eqref{N-boundary}. How to get rid of this
condition has its independent interest.
\end{remark}

We are now in a position to give the proof of Theorem \ref{th-2}.

\ms

{\em Proof of Theorem \ref{th-2}}.  For simplicity,  we always
assume that $T_0 = T$ in this section. If $y(\cdot,T) = 0$ in $G$,
$P$-a.s., then the inequality \eqref{th-2-eq} holds. Now we only
consider the case that $y(\cdot,T) \neq 0$ in $G$, $P$-a.s. In this
case, due to Lemma \ref{inv th1}, we know that $y(\cd,t)\neq 0$  in
$G$, $P$-a.s.

 Note that $B_r\subset G_0 \subset\subset G$ is the ball
centered at $x_0$ with radius $r$. As in
the proof for Theorem \ref{th-1}, we begin
with calculating
\begin{equation*}
  \int_G \|\nabla \left[y(x,T) \vartheta(2\lambda, x,x_0)\right]\|^2 dx \geq 0.
\end{equation*}
Via some straightforward calculations,
similar to the process of getting
\eqref{N15.2}, we arrive at that
\begin{equation}
  \label{N-convex-1}
\begin{array}
    {ll}
    &\quad \ds \[1- \frac{8\lambda}{r^2}\(\lambda\widetilde N(T)+\frac{n}{2}\)\]
    \mE \int_G |x - x_0|^2 y^2(x,T) \vartheta(\lambda, x,x_0)dx\\
    \ns   & \leq \ds 8\lambda \(\lambda \widetilde N(T)+\frac{n}{2}\) \mE \int_{B_r} y^2(x,T) \vartheta(\lambda, x,x_0)dx.
  \end{array}
\end{equation}

In what follows, we turn to estimate $\lambda \widetilde N(T)$.

From the estimate for $\widetilde N$ in Lemma \ref{lemma 4} and
Gronwall's inequality, it follows that
\begin{equation}
  \label{N-convex-2}
  \begin{array}
  {ll}
  &\ds \quad \lambda \widetilde N(T)\\
  \ns  & \ds \leq \exp\big({(T -t)\norm{b}^2_{G}}\big)(T -t+ \lambda)\widetilde N(t) \\
   \ns \ms & \quad \ds + \exp\big({(T -t)\norm{b}^2_{G}}\big)(T -t + \lambda)(T -t) \big( \norm{a}^2_{G}
    + 2\norm{b}^2_{G}\big)\\
    \ns \ms & \ds \leq \exp\({T\norm{b}^2_{G}}\)(T + \lambda)\widetilde N(t)
     + \exp\big({T\norm{b}^2_{G}}\big)(T + \lambda)T  \big( \norm{a}^2_{G}
    + 2\norm{b}^2_{G}\big).
    \end{array}
\end{equation}

Integrating \eqref{N-convex-2} with time variable $t$ over $[0, \frac
T2]$, we have that
\begin{equation*}
  \begin{array}
    {ll}
\ds \q\frac T 2 \lambda \widetilde N(T)  \\
\ns\ds\ds \leq
(T+\lambda)\exp\big({T\norm{b}^2_{G}}\big)\int_0^{\frac
  T2} \widetilde N(\tau)d\tau  + (T+\lambda)\frac{T^2}{2}\exp\big({T\norm{b}^2_{G}}\big)\big(\norm{a}^2_{G} +
    2\norm{b}^2_{G}\big).
  \end{array}
\end{equation*}
Turning to the result of Lemma \ref{lemma
3}, say \eqref{H-convex}, it is easy to
show that
\begin{equation*}
  \int_{0}^{\frac T2} \widetilde N(\tau)d\tau \leq \ln
  \frac{\widetilde H(0)}{\widetilde H(\frac T2)} + T \big(\norm{a}_{G}
   + 2\norm{b}^2_{G}\big).
\end{equation*}
Thus, we obtain that
\begin{equation}
  \label{N-convex-3}
  \begin{array}
    {ll}
    \ds \hspace{-0.4cm}  \lambda \widetilde N(T)\3n & \ds \leq
\frac{2(T+\lambda)}{T}
\exp\big({T\norm{b}^2_{G}}\big)\[\ln\frac{\widetilde
H(0)}{\widetilde H(\frac T2)} + \frac T2
\big(2 \norm{a}_{G} + T \norm{a}^2_{G} + (1
+ 2T)\norm{b}^2_{G}\big)\].
  \end{array}
\end{equation}

Denote $m = \max_{x \in \overline{G}}|x - x_0|^2$. By some elementary
calculations, we find that
\begin{equation}
  \label{N-convex-7}
  \begin{array}
  {ll}
\ds \frac{\widetilde H(0)}{\widetilde
H(\frac T2)} \3n& \ds = \left(\frac{\frac
T2 + \lambda}{T + \lambda}\right)^{\frac n
2}
  \frac{\ds\mE\int_G y^2(x, 0)\vartheta(T+\lambda, x,x_0)dx}
  {\ds\mE\int_G y^2\(x, \frac T2\)\vartheta(T/2 +\lambda, x,x_0)dx}\\
  \ns & \ds \leq \exp\(\frac m{4(\lambda + \frac T2)}\) \frac{\ds\mE\int_G y^2(x, 0)dx}
  {\ds\mE\int_G y^2\(x, \frac T2\)dx}\\
  \ns & \leq \ds  \exp\left(\frac m{2T}\right) \frac{\ds\mE\int_G y^2(x, 0)dx}
  {\ds\mE\int_G y^2\(x, \frac T2\)dx}.
  \end{array}
\end{equation}

Note that $y$ solves the equation
\eqref{system1}. By It\^o's formula, we have
that
\begin{equation*}
  \begin{array}
  {ll}
    & \quad\ds \mE\int_G y^2(x,t)dx \\
    \ns  & \ds  = \mE\int_G y^2(x,s)dx \!+ \mE\int_s^t \int_G
    (2ay^2 + b^ 2y^2 - 2|\nabla y|^2)dxd\tau\\
    \ns & \ds \leq \mE\int_G y^2(x,s)dx + \mE\int_s^t \int_G
    (2ay^2 + b^ 2y^2)dxd\tau\\
    \ns & \ds \leq  \mE\int_G y^2(x,s)dx \!+\! \(2\norm{a}_{G}
    +
     \norm{b}^2_{G}\)\mE\int_s^t \int_G
    y^2dxd\tau
  \end{array}
\end{equation*}
holds for $0\leq s \leq t \leq T$. Due to Gronwall's inequality, it
follows that
\begin{equation}
\label{N-convex-8}
  \mE \int_G y^2(x,T)dx \leq \exp\(T\(\norm{a}_{G}
    +
     \frac {\norm{b}^2_{G} } 2\)\)\mE
     \int_G y^2\(x,\frac T2\)dx.
\end{equation}

According to \eqref{N-convex-7} and \eqref{N-convex-8}, we obtain
that
\begin{equation}
  \label{N-convex-9}
  \frac{\widetilde H(0)}{\widetilde H(\frac T2)} \leq
  \exp\(\frac m{2T} + T\(\norm{a}_{G}
    +
     \frac {\norm{b}^2_{G} } {2}\)\)\frac{\ds \mE\int_G y^2 (x, 0) dx}{\ds\mE \int_G
     y^2(x,T)dx}.
\end{equation}

Making use of \eqref{N-convex-3} and \eqref{N-convex-9}, it is easy
to show that
\begin{equation}
  \label{N-convex-10}
  \begin{array}
    {ll}
    \ds   \lambda \widetilde N(T)\3n & \ds \leq
(T+\lambda)\exp({T\norm{b}^2_{G}})\[\frac 2T\ln\frac{\ds \mE \int_G
y^2(x,0)dx}{\ds \mE \int_G y^2(x,T)dx} + \frac{m}{T^2} \\
\ns & \ds \quad\!\! +  \big(4 \norm{a}_{G}
+ T \norm{a}^2_{G} + 2(1 +
T)\norm{b}^2_{G}\big)\].
  \end{array}
\end{equation}

Fix $\lambda \in (0, 1]$. Denote
\begin{equation*}
\begin{array}
    {ll}
    \ds  \cD \3n & \ds =(T+1)\exp\big({T\norm{b}^2_{G}}\big)\[\frac 2T\ln\frac{\ds \mE \int_G
y^2(x,0)dx}{\ds \mE \int_G y^2(x,T)dx} \\
\ns & \ds \quad\!\! + \frac{m}{T^2} +
\big(4 \norm{a}_{G} + T \norm{a}^2_{G} +
2(1 + T)\norm{b}^2_{G}\big)\] + \frac n 2.
  \end{array}
\end{equation*}
We  find that
\begin{equation*}
  \lambda \widetilde N(T) + \frac n2 \leq \cD.
\end{equation*}
This together with \eqref{N-convex-1}
implies
\begin{equation*}
\begin{array}
{ll}
  & \ds  \(1- \frac{8\lambda}{r^2}\cD \)\mE \int_G |x - x_0|^2 y^2(x,T)
  \vartheta(\lambda,x,x_0)dx\leq 8\lambda \cD \mE \int_{B_r} y^2(x,T) \vartheta (\lambda,x,x_0)dx.
  \end{array}
\end{equation*}
Letting $\widetilde \lambda = \frac{r^2}{16\cD}$, we arrive at that
\begin{equation}
  \label{N-convex-11}
\mE \int_G |x - x_0|^2 y^2(x,T) \vartheta (\widetilde\lambda,x,x_0)dx
\leq r^2 \mE \int_{B_r} y^2(x,T) \vartheta
(\widetilde\lambda,x,x_0)dx.
\end{equation}
Due to \eqref{N-convex-11}, an elementary computation gives that
\begin{equation}
  \label{N-convex-12}
  \begin{array}
    {ll}
\ds \q\mE\int_G y^2(x,T) \vartheta(\widetilde\lambda,x,x_0)dx  \\
  \ns\ds = \mE\int_{G\setminus B_r} y^2(x,T) \vartheta (\widetilde\lambda,x,x_0)dx
  + \mE\int_{ B_r} y^2(x,T) \vartheta (\widetilde\lambda,x,x_0)dx\\
  \ns   \ds \leq \frac 1{r^2}\mE\int_G |x - x_0|^2 y^2(x,T) \vartheta(\widetilde\lambda,x,x_0)dx
  + \mE\int_{ B_r} y^2(x,T)\vartheta(\widetilde\lambda,x,x_0)dx\\
  \ns  \ds \leq 2 \mE\int_{ B_r} y^2(x,T) \vartheta(\widetilde\lambda,x,x_0)dx.
  \end{array}
\end{equation}

Noting the definition of $m$ and the choice of $\widetilde \lambda$,
we follows from \eqref{N-convex-12} that
\begin{equation}
\label{N-convex-13}
  \mE \int_G y^2(x,T)dx \leq 2 \exp\(\frac{4m\cD}{r^2}\)\mE \int_{B_r}
  y^2 (x,T)dx.
\end{equation}

Put
\begin{equation*}
\begin{array}{ll}\ds
  \cJ \3n&= \cD - \frac{2(T+1)}{T} \exp({T\norm{b}^2_{G}})\ln\frac{\ds \mE \int_G y^2(x,0)dx}{\ds \mE
\int_G y^2(x,T)dx}\\
\ns&\ds =
(T+1)\exp\big({T\norm{b}^2_{G}}\big)\[
\frac{m}{T^2} +  \big(4 \norm{a}_{G} + T
\norm{a}^2_{G} + 2(1 +
T)\norm{b}^2_{G}\big)\] + \frac n 2.
\end{array}
\end{equation*}
Plugging $\cJ$ into \eqref{N-convex-13},
one has that
\begin{equation*}
\begin{array}
{ll}
 & \ds  \q\mE\int_G y^2(x,T)dx \\
 \ns & \ds \leq 2 \exp\left(\frac{4m\cJ}{r^2}\right)\left(\frac{ \ds \mE\int_G y^2(x,0)dx}{ \ds \mE\int_G
y^2(x,T)dx}\right)^{\theta}
\mE\int_{B_r}y^2(x,T)dx,
\end{array}
\end{equation*}
where $\theta = \frac{8m(T+1)\exp({T\norm{b}^2_{G}})}{r^2T}$.

 After some elementary calculations,
one finds that
\begin{equation}
  \label{N-convex-14}
  \mE\int_G y^2(x,T)dx \leq 2^{\delta} \exp(\beta)\left(\mE \int_G y^2(x,0)dx
  \right)^{1 - \delta} \left(\mE \int_{B_r} y^2(x,T)dx\right)^{\delta}.
\end{equation}
In \eqref{N-convex-14}, the simple notations
$\delta$ and $\beta$ are given as follows:
\begin{equation*}
  \begin{cases}
    \ds \delta = \frac{r^2 T}{r^2 T  + 8m (T + 1)\exp(T\norm{b}^2_{G})},\\
\ns \ds \beta = \frac{4m T \cJ}{r^2 T  + 8m
(T + 1)\exp(T\norm{b}^2_{G})}.
  \end{cases}
\end{equation*}
Noting that $B_r \subset G_0$, we then complete the proof of Theorem
\ref{th-2}.
\endpf

\section{Proof of Theorem \ref{observability theorem}}

\label{proof of observability}

This section is devoted to the proof of
Theorem \ref{observability theorem}. To
begin with, we first introduce the
following lemma, whose proof can be found
in \cite{phung2}.

\begin{lemma}
\label{ob-lemma}
  Let $E \subset (0, T)$ be a measurable set of positive measure. Let
  $t_0 $ be a density point of $E$. Then for each $z > 1$, there
  exists an $t_1 \in (t_0, T)$ such that the sequence
  $\{t_m\}_{m=1}^{\infty}$, given by
  \begin{equation}
    t_{m+1} = t_0 + \frac 1{z^m} (t_1  - t_0),
  \end{equation}
  satisfies
  \begin{equation}
    t_m  - t_{m+1}\leq 3|E\cap (t_{m+1}, t_m)|.
  \end{equation}
\end{lemma}

\ms \ms

{\em Proof of Theorem \ref{observability
theorem}}:  Let $t_0$ be a density point of
$E$. Let $\{t_m\}_{m=1}^{\infty}$ be a
sequence provided by Lemma \ref{ob-lemma}.
From \eqref{N-convex-14}, it follows that
for any $t\in [0,T]$,
\begin{equation}
  \label{ob-1}
  \mE\int_G y^2(x,t)dx \leq \frac 2{\varepsilon^{\gamma}}\exp(\Theta) \mE\int_{B_r}y^2(x,t)dx + \varepsilon \mE \int_G
  y^2(x,0)dx,
\end{equation}
where
\begin{equation*}
\Theta = \max_{t\in
[0,T]}\frac{t\cJ}{2(t+1)\exp\big(t\norm{b}^2_{G}\big)},\;\;
\gamma =\max_{t\in [0,T]}\frac{8m (t +
1)\exp\big(t\norm{b}^2_{G}\big)}{r^2 t  +
8m (t + 1)\exp\big(t\norm{b}^2_{G}\big)}.
\end{equation*}

For the solution to \eqref{system1}, by standard argument, it is easy
to verify the following energy estimate:
\begin{equation}
\label{energy estimate} \mE\int_G
y^2(x,t)dx \leq \exp\big(C(a,b)t\big)\mE
\int_G y^2(x,0)dx,
\end{equation}
where
\begin{equation*}
C(a,b, t) = (2\norm{a}^2_{G} + \norm{b}^2_{G})t.
\end{equation*}

For $t \in (t_{m+1}, t_m]$, combining
\eqref{ob-1} and \eqref{energy estimate},
one finds that
\begin{equation*}
  \begin{array}
    {ll}
    \ns \ds\exp (-C(a,b, T))\mE\int_G y^2(x,t_m)dx \3n& \leq \ds \mE\int_G y^2(x,t)dx \\
    \ns & \ds \leq \frac 2{\varepsilon^{\gamma}}\exp(\Theta) \mE \int_{B_r}y^2(x,t)dx  + \varepsilon \mE\int_G y^2(x,t_{m+1})dx.
  \end{array}
\end{equation*}
As a direct result of this inequality, we get that
\begin{equation*}
\varepsilon^{\gamma}\exp (-C(a,b, T))\mE\int_G y^2(x,t_m)dx -
\varepsilon^{\gamma + 1}\mE\int_G y^2(x,t_{m+1})dx \leq 2
\exp(\Theta) \mE \int_{B_r}y^2(x,t)dx.
\end{equation*}
Integrating over $E\cap (t_{m+1}, t_m)$ shows that
\begin{equation}
\begin{array}
{ll}
  & \ds |E\cap (t_{m+1}, t_m)| \[\varepsilon^{\gamma}\exp (-C(a,b, T))\mE\int_G
  y^2(x,t_m)dx - \varepsilon^{\gamma + 1}\mE\int_G
  y^2(x,t_{m+1})dx\]\\
  \ns & \ds \leq  2 \exp(\Theta)\mE\int_{E\cap (t_{m+1}, t_m)}\int_{B_r}y^2(x,t)dx.
  \end{array}
\end{equation}

  Letting $\alpha_m = \varepsilon ^{\gamma}|E\cap (t_{m+1}, t_m)| , \sigma_m  = \varepsilon^{\gamma + 1}|E\cap (t_{m+1},
  t_m)|$, then it is clear that
  \begin{equation*}
  \label{ob-2}
\begin{array}
{ll}
  & \ds \alpha_m  \exp (-C(a,b, T))\mE\int_G
  y^2(x,t_m)dx - \sigma_m \mE\int_G
  y^2(x,t_{m+1})dx\\
  \ns & \ds \leq  2 \exp(\Theta)\mE\int_{E\cap (t_{m+1}, t_m)}\int_{B_r}y^2(x,t)dx.
  \end{array}
\end{equation*}

We choose a sequence of $\{\varepsilon_m\}_{m=1}^{\infty}$ in the
following way: Let $\sigma_m = \alpha_{m+1} \exp(-C(a,b, T))$ such
that
$$
\frac{\varepsilon_m^{\gamma +
1}}{\varepsilon_{m+1}^{\gamma}}\exp(C(a,b,T))
= \frac{|E\cap (t_{m+2}, t_{m+1})|}{|E\cap
(t_{m+1}, t_m)|}.
$$
By means of the properties of $\{t_m\}_{m=1}^{\infty}$ mentioned in
Lemma \ref{ob-lemma}, this implies that for any $m\in\dbN$,
\begin{equation}\label{3.26-eq1}
\begin{array}{ll}\ds
\varepsilon_{m+1}^{\gamma}\3n&\ds =
\varepsilon_m^{\gamma +
1}\exp(C(a,b,T))\frac{|E\cap (t_{m+1},
t_m)|}{|E\cap (t_{m+2}, t_{m+1})|}\\
\ns&\ds \leq \varepsilon_m^{\gamma +
1}\exp(C(a,b,T)) \frac{3|(t_{m+1},
t_m)|}{|(t_{m+2}, t_{m+1})|}\\
\ns&\ds \leq 3z\varepsilon_m^{\gamma +
1}\exp(C(a,b,T)).
\end{array}
\end{equation}
Let us choose $\e_1=\frac{1}{3z\exp(C(a,b, T))}$. Then from
\eqref{3.26-eq1}, we see that $\e_2^{\g}\leq \e_1^{\g}$. Hence, we
know that $\varepsilon_2\leq\frac{1}{3z\exp(C(a,b, T))}$. This,
together with \eqref{3.26-eq1}, implies that $\varepsilon_3
\leq\frac{1}{3z\exp(C(a,b))}$. By induction, we find that
\begin{equation*}
   \e_m\leq
  \frac{1}{3z\exp(C(a,b,T))} \mbox{ for any }m\in\dbN.
\end{equation*}
Consequently, following \eqref{ob-2}, one finds that
\begin{equation*}
\begin{array}
{ll}
  & \ds \sum_{m=1}^n \left(\alpha_m  \exp (-C(a,b,T))\mE\int_G
  y^2(x,t_m)dx - \sigma_m \mE\int_G
  y^2(x,t_{m+1})dx\right)\\
  \ns & \ds \leq  2 \exp(\Theta)\sum_{m=1}^n \mE\int_{E\cap (t_{m+1},
  t_m)}\int_{B_r}y^2(x,t)dx,
  \end{array}
\end{equation*}
which implies that
\begin{equation}
\label{ob-3}
\begin{array}
{ll}
  & \ds \alpha_1  \exp (-C(a,b,T))\mE\int_G
  y^2(x,t_1)dx - \sigma_n \mE\int_G
  y^2(x,t_{n+1})dx\\
  \ns & \ds \leq  2 \exp(\Theta)
  \mE\int_{E}\int_{B_r}y^2(x,t)dx.
  \end{array}
\end{equation}

It is noted that
\begin{equation*}
  \lim_{n\to \infty} \sigma_n  = \lim_{n\to \infty}  \varepsilon_n^{\gamma +
  1}|E\cap (t_{n+1}, t_n)|  = 0
\end{equation*}
following the construction of
$\{t_n\}_{n=1}^{\infty}$ and that $0 <
\varepsilon _n < \frac{1}{3z\exp(C(a,b,
T))}$. Therefore, in light of \eqref{ob-3},
we immediately arrive at
\begin{equation*}
\label{ob-4}
\begin{array}
{ll}
  & \ds\mE\int_G
  y^2(x,t_1)dx \leq  2  \alpha_1^{-1}  \exp (C(a,b,T))\exp(\Theta)
  \mE\int_{E}\int_{B_r}y^2(x,t)dx.
  \end{array}
\end{equation*}
This inequality together with the energy estimate \eqref{energy
estimate} shows that

\begin{equation*}
\label{ob-5}
\begin{array}
{ll}
  & \ds\mE\int_G
  y^2(x,T)dx \leq  2  \alpha_1^{-1}  \exp (2C(a,b,T))\exp(\Theta)
  \mE\int_{E}\int_{B_r}y^2(x,t)dx.
  \end{array}
\end{equation*}
Letting $C =  2  \alpha_1^{-1}  \exp (2C(a,b,T))\exp(\Theta)$ and
noting that $B_r \subset G_0$, we complete the desired result. \endpf

\section{Proof of Theorem \ref{th-3} and Theorem \ref{null control}}
\label{proof of th3}

{\it Proof of Theorem \ref{th-3}}\,: Due to the linearity of
\eqref{csystem2}, it suffices to show that its attainable set $A_T$
at time $T$ with final datum $z(T)=0$ is dense in $L^2(G)$. Let us
prove this by the contradiction argument.

If $A_T$ is not dense in $L^2(G)$, then there exists an $\eta\in
L^2(G)$ such that $\eta \neq 0$ and
$$ \int_G z(0)\eta dx=0 \mbox{ for any
}z(0)\in A_T.
$$
Let us consider the following equation
\begin{equation}\label{system2}
\left\{
\begin{array}{ll}\ds
d\tilde y - \D \tilde y dt = -a_1\tilde y dt -
b_1\tilde ydB(t) &\mbox{ in } G\times (0,T),\\
\ns\ds \tilde y=0 &\mbox{ on } \partial G\times (0,T),\\
\ns\ds \tilde y(0)=\eta &\mbox{ in }G.
\end{array}
\right.
\end{equation}
It is clear that the solution $\tilde y$ to \eqref{system2} belongs
to $L^2_\cW(\O;C([0,T];L^2(G)))\cap L^2_\cW(0,T;H_0^1(G))$.

From It\^{o}'s formula, we find that
\begin{equation}\label{6.9-eq1}
\begin{array}{ll}
\ds \q\mE\int_G \tilde y(T)z(T)dx-\mE\int_G \tilde y_0z(0)dx\\
\ns\ds = \mE \int_0^T\int_G \[\tilde y \big(-\D z + a_1z + b_1 Z + h
+ \chi_{E_1}\chi_{G_0}f \big) + z\big( \D y - a_1 \tilde y\big) - b_1
\tilde y Z
\]dxdt
\\
\ns \ds=\mE\int_0^T\int_{G}\tilde yh\, dxdt +
\mE\int_0^T\int_{G}\tilde y\chi_{E_1}\chi_{G_0}f\, dxdt.
 \end{array}
 \end{equation}
Let $z(T)=0$, we obtain that
\begin{equation}\label{app1}
0=\mE\int_G z(0)\eta dx = \mE\int_0^T\int_{G}\tilde yh\, dxdt
+\mE\int_0^T\int_{G}\tilde y\chi_{E_1}\chi_{G_0}f dxdt.
\end{equation}
Hence,
$$
\mE\int_0^T\int_{G}\tilde yh\, dxdt +\mE\int_0^T\int_{G}\tilde
y\chi_{E_1}\chi_{G_0}f dxdt =0
$$
for any $f\in L^2_{\cW}(0,T;L^2(G))$. Therefore we get
$$
\tilde y=0 \mbox{ in } G_0\times E_1,\; P\mbox{-a.s.}
$$
Then, from Theorem \ref{th-1},  we arrive at $\eta=0$, a
contradiction.
\endpf

{\it Proof of Theorem \ref{null control}}\,: Consider the following
equation:
\begin{equation}\label{system3}
\left\{
\begin{array}{ll}\ds
d\hat y - \D \hat y dt = -a_1\hat y dt -
b_1\hat ydB(t) &\mbox{ in } G\times (0,T),\\
\ns\ds \hat y=0 &\mbox{ on } \partial G\times (0,T).
\end{array}
\right.
\end{equation}
We introduce a linear subspace of $L_{\cW}^2(0,T;L^2(G_0))$:
$$
\cX\=\big\{\hat y|_{G_0\times E_1}\;:\;\hat y\hb{ solves
 equation
 }\eqref{system3}\mbox{ with some initial datum } \hat y_0\in
L^2(G)\big\},
$$
and define a linear functional $\cL$ on
$\cX$ as follows:
 $$
 \cL(\hat y|_{G_0\times
E_1})=-\mathbb{E}\int_G y(T)z_Tdx.
 $$
By Theorem \ref{observability theorem}, we see that
$$
\begin{array}{ll}\ds
|\cL(\hat y|_{G_0\times E_1})| \3n&\ds\leq
\norm{y(T)}_{L^2(\O,\cF_T,P;L^2(G))}

\norm{z_T}_{L^2(\O,\cF_T,P;L^2(G))}\\
\ns&\ds  \leq
\cC\norm{z_T}_{L^2(\O,\cF_T,P;L^2(G))}\(\mE\int_{0}^T\int_{G_0}|\hat
\chi_{E_1}y|^2dxdt\)^{\frac{1}{2}}.
\end{array}
$$
Therefore, $\cL$ is a bounded linear functional on $\cX$. By
Hahn-Banach Theorem, $\cL$ can be extended to a bounded linear
functional with the same norm on $L_{\cF}^2(0,T;L^2(G_0))$. For
simplicity, we use the same notation for this extension. Now, Riesz
Representation Theorem allows us to find a  random field $f\in
L_{\cF}^2(0,T;L^2(G_0))$ so that
\begin{equation}\label{7.2}
\mE\int_0^T\int_{G}\hat y\chi_{E_1}\chi_{G_0}f\, dxdt
=\mathbb{E}\int_G y(T)z_Tdx.
\end{equation}

We claim that this $f$ is the control we need. In fact, for any
$z_T\in L^2_{\cW_T}(\Omega;L^2(G))$, for the solution $\hat y$ of
equation \eqref{system3} and the solution $(z,Z)$ of equation
\eqref{csystem2}, by It\^o's formula, we get that
\begin{equation}\label{6.9-eq11}
\begin{array}{ll}
\ds \q\mE\int_G \hat y(T)z(T)dx-\mE\int_G \hat y_0z(0)dx\\
\ns\ds = \mE \int_0^T\int_G \[\hat y \big(-\D z + a_1z + b_1 Z  +
\chi_{E_1}\chi_{G_0}f \big) + z\big( \D y - a_1 \hat y\big) - b_1
\hat y Z
\]dxdt
\\
\ns \ds=  \mE\int_0^T\int_{G}\hat y\chi_{E_1}\chi_{G_0}f\, dxdt.
\end{array}
\end{equation}
Combining \eqref{7.2} and \eqref{6.9-eq11}, we find that
 $$
 \mE\int_G \hat y_0z(0)dx=0.
 $$
Since $\hat y_0$ can be chosen arbitrarily, we know that $z(0)=0$ in
$G$.
\endpf

\section*{Acknowledgement}

Qi L\"{u} would like to thank the Laboratoire Jacques-Louis Lions
(LJLL) at Universit\'{e} Pierre et Marie Curie for its hospitality.
Part of this work was carried out while he was visiting the LJLL.




\end{document}